\ifx\LocalElevenPtArticleFormatLoaded\undefined
  \documentclass[a4paper,11pt]{article}     
  \usepackage{array}                        
  \usepackage{theorem}                      
  \usepackage{amsmath,amscd,amssymb}        
  \usepackage{latexsym}                     
  \usepackage[german,english]{babel}        
  \usepackage[applemac]{inputenc}
\usepackage{epsfig}
\usepackage{xypic}
\fi

\newtheorem{theorem}{Theorem}[section]
\usepackage{theorem}

\newtheorem{lemma}[theorem]{Lemma}

\newtheorem{proposition}[theorem]{Proposition}

\theorembodyfont{\rm}
\newtheorem{remark}{Remark}

\newcommand{\CC}{{\mathbb{C}}}

\newcommand{\HH}{{\mathbb{H}}}

\newcommand{\PP}{{\mathbb{P}}}
\newcommand{\QQ}{{\mathbb{Q}}}
\newcommand{\RR}{{\mathbb{R}}}
\newcommand{\ZZ}{{\mathbb{Z}}}

\newcommand{\NN}{{\mathbb{N}}}

\newenvironment{Proof}{\begin{ProofwCaption}{Proof}}{\end{ProofwCaption}}
\newenvironment{Proof*}[1]{\begin{ProofwCaption}{{#1}}}{\end{ProofwCaption}}
\newenvironment{ProofwCaption}[1]%
  {\addvspace\theorempreskipamount \noindent{\it #1.}\rm}%
  {\qed \par \addvspace\theorempostskipamount}
\newcommand{\qedsymbol}{\mbox{$\Box$}}
\newcommand{\qed}{\quad\qedsymbol}
\begin{document}
\begin{center}
{\bf\Large Igusa's Modular Form and the
Classification of Siegel Modular Threefolds}\\[3mm]
{\large Klaus Hulek}
\end{center}

\section{Introduction}
For an integer $d\ge 1$ let
$$
E_d=\begin{pmatrix}
1 & 0\\ 0 & d\end{pmatrix}, \quad
\Lambda_d=
\begin{pmatrix}
0    & E_d\\
-E_d & 0
\end{pmatrix}.
$$
We consider the symplectic group
$$
\tilde{\Gamma}_{1,d}:=\operatorname{Sp}(\Lambda_d,\ZZ).
$$
For $d=1$ this is the usual integer symplectic group
$\operatorname{Sp}(4,\ZZ)$. The group
$\tilde{\Gamma}_{1,d}$ operates on the {\em Siegel space} of genus 2
$$
\HH_2=\left\{
\tau=
\begin{pmatrix}
\tau_1 & \tau_2\\
\tau_2 & \tau_3
\end{pmatrix}
\in \operatorname{Mat} (2\times 2,\CC); \mbox{Im }\tau>0
\right\}
$$
by
$$
\tilde{M}=
\begin{pmatrix}
\tilde{A} & \tilde{B}\\
\tilde{C} & \tilde{D}
\end{pmatrix}
:\tau\mapsto (\tilde{A} \tau+\tilde{B} E_d) (\tilde{C} \tau+\tilde{D}
E_d)^{-1} E_d.
$$
The quotient
$$
{\cal{A}}_{1, d}=\tilde{\Gamma}_{1,d}\backslash \HH_2
$$
is the moduli space of $(1,d)$-polarized abelian surfaces. Alternatively we
can consider the following subgroup of the
usual rational symplectic group $\operatorname{Sp}(4,\QQ)$. Let
$$
R_d=\mbox{diag} (1,1,1,d)
$$
and set
$$
\Gamma_{1,d}:=R^{-1}_d\tilde{\Gamma}_{1,d} R_d\subset \operatorname{Sp}
(4,\QQ).
$$
Then $\Gamma_{1,d}$ acts in the usual way on $\HH_2$ by
$$
M=\begin{pmatrix}
A & B\\
C & D
\end{pmatrix}:\tau\mapsto (A\tau+B) (C\tau+D)^{-1}
$$
and
$$
{\cal A}_{1,d}=\tilde{\Gamma}_{1,d}\backslash\HH_2=\Gamma_{1,d}\backslash
\HH_2.
$$
Let $L=\ZZ^4$ be the lattice on which $\Lambda_d$ defines a symplectic form
and let $L^{\vee}$ be the dual lattice of
$L$. We consider the following subgroups of
$\tilde{\Gamma}_{1,d}$ defined by
$$
\begin{array}{lcl}
\tilde{\Gamma}_{1,d}^{\operatorname{lev}} & :=
&\left\{M\in\tilde{\Gamma}_{1,d};\quad M|_{L^{\vee}/L}=\mbox{id}\right\}\\
[2mm]
\tilde{\Gamma}_{1,d}(n)                   & :=
&\left\{M\in\tilde{\Gamma}_{1,d};\quad M\equiv{\bf 1} \mbox{ mod
}n\right\}\quad (n\ge 1)\\[2mm]
\tilde{\Gamma}_{1,d}^{\operatorname{lev}}(n) & :=
&\tilde{\Gamma}_{1,d}^{\operatorname{lev}}
\cap \tilde{\Gamma}_{1,d}(n).
\end{array}
$$
This gives rise to subgroups of $\operatorname{Sp} (4,\QQ)$:
$$
\begin{array}{lcl}
{\Gamma}_{1,d}^{\operatorname{lev}}    & :=
&R^{-1}_d\tilde{\Gamma}_{1,d}^{\operatorname{lev}} R_d\\[2mm]
{\Gamma}_{1,d}(n)                      & :=
&R_d^{-1}\tilde{\Gamma}_{1,d}(n)R_d\\[2mm]
{\Gamma}_{1,d}^{\operatorname{lev}}(n) & :=
&R_d^{-1}\tilde{\Gamma}_{1,d}^{\operatorname{lev}}(n)R_d=
{\Gamma}_{1,d}^{\operatorname{lev}}\cap \Gamma_{1,d}(n),
\end{array}
$$
resp. to the moduli spaces
$$
\begin{array}{lclcl}
{\cal
A}_{1,d}^{\operatorname{lev}}&=&\tilde{\Gamma}_{1,d}^{\operatorname{lev}}
\backslash\HH_2
&=&{\Gamma}_{1,d}^{\operatorname{lev}}\backslash\HH_2\\
{\cal A}_{1,d}(n)&=&\tilde{\Gamma}_{1,d}(n) \backslash\HH_2
&=&{\Gamma}_{1,d}(n)\backslash\HH_2\\
{\cal
A}_{1,d}^{\operatorname{lev}}(n)&=&\tilde{\Gamma}_{1,d}^{\operatorname{lev}}(n)
\backslash\HH_2
&=&{\Gamma}_{1,d}^{\operatorname{lev}}(n)\backslash\HH_2 .
\end{array}
$$
The geometric meaning of these moduli spaces is the following:\\[3mm]
\begin{tabular}{lr}
{ } &${\cal A}_{1,d}^{\operatorname{lev}}=\{(A, H, \alpha); (A,H)$ is a
$(1,d)$-polarized abelian surface,\\
{ } &$\alpha$ is a canonical level-structure\}.
\end{tabular}
\\[2mm]
Here a {\em canonical} level-structure is a symplectic basis of the kernel
of the map $\lambda_H:A\rightarrow
\hat{A}=\mbox{Pic}^{0} A$. (Note that this kernel is (non-canonically)
isomorphic to $\ZZ/d\times \ZZ/d$.)
Similarly\\[3mm]
\begin{tabular}{lr}
{ } &${\cal A}_{1,d}(n)=\{(A, H, \beta); (A,H)$ is a $(1,d)$-polarized
abelian surface,\\
{ } &$\beta$ is a full level-n structure\}.
\end{tabular}\\[2mm]
Here a {\em full level}-n {\em structure} is a symplectic basis of the
group $A^{(n)}$
 of $n$-torsion points of $A$. Finally\\[3mm]
\begin{tabular}{lr}
{ } &${\cal A}^{\operatorname{lev}}_{1,d}(n)=\{(A, H,\alpha, \beta); (A,H)$
is a $(1,d)$-polarized abelian surface,\\
{ } &$\alpha$ is a canonical level structure, $\beta$ is a full level-$n$
structure\}.
\end{tabular}
\\[2mm]
Note that
$$
\tilde{\Gamma}_{1,d}^{\operatorname{lev}}/
 \tilde{\Gamma}_{1,d}\cong {\Gamma}_{1,d}^{\operatorname{lev}}/
{\Gamma}_{1,d}\cong\operatorname{SL}(2,\ZZ/d)
$$
and that we have, therefore, a Galois covering ${\cal
A}_{1,d}^{\operatorname{lev}}\rightarrow {\cal A}_{1,d}$
 with Galois group $\mbox{PSL}(2,\ZZ/d)$.

The aim of this short note is to prove two results about the classification
of these Siegel modular varieties.

\begin{theorem}
${\cal A}_{1,d}(n)$ is of general type if $(d,n)=1$ and $n\ge 4$.
\end{theorem}

Since ${\cal A}_{1,1}(n)$ is rational for $n\le 3$ this is the best result
which one can hope for if one considers all
$d$ simultaneously. The space ${\cal A}_{1,3}(2)$ has a Calabi-Yau model
(\cite{BN}, \cite{GH}) and hence Kodaira
dimension $0$, whereas  ${\cal A}_{1,3}(3)$ is of general type
(\cite[Theorem 3.1]{GH}). For prime numbers $p$
Sankaran \cite{S} has proved that ${\cal A}_{1,p}$ is of general type for
$p\ge 173$. A similar result for ${\cal
A}_{1,d}$, where $d$ is not necessarily prime, is, as for as I know, not
known. Borisov \cite{Bo} has shown that, up to
conjugation, there are only finitely many subgroups $\Gamma$ of
$\operatorname{Sp}(4,\ZZ)$ such that ${\cal
A}(\Gamma)=\Gamma\backslash\HH_2$ is not of general type. Recall however,
that the groups $\Gamma_{1,d}(n)$ are not
subgroups of
$\operatorname{Sp}(4,\ZZ)$ unless $d$ divides $n$ and that, in general,
they are also not conjugate to subgroups of
$\operatorname{Sp}(4,\ZZ)$. An essential ingredient in the proof of the
theorem is Igusa's modular form ${\Delta}_{10}$.

The above theorem is a result about the birational classification of these
varieties. If one wants to ask more
precise questions, such as whether $K$ is ample, then one has to specify
the compactification with which one wants to
work.

\begin{theorem}
The Voronoi compactification $\left(
{\cal A}_{1,p}^{\operatorname{lev}}
(n)\right)^{\ast}$ for a prime number $p$ with $(p,n)=1$ is smooth and has
ample canonical bundle (i.e. is a canonical
model in the sense of Mori theory)
 if and only if $n\ge 5$.
\end{theorem}

Here a few words are in order: By Voronoi compactification we mean the
compactification defined by the second Voronoi
decomposition. We choose this compactification because Alexeev
\cite{A1} has shown that it appears naturally when one wants to construct a
toroidal compactification which represents
a geometrically meaningful functor. The addition of a canonical level
structure has two reasons: The spaces $\left({\cal
A}_{1,p}(n)\right)^{\ast}$ have non-canonical singularities for infinitely
many $p$ and $n$. These singularities come
from the toroidal construction, not from fixed points of the group which is
neat for $n\ge 3$. Moreover, it is necessary
to introduce at least some kind of level structure to obtain a functorial
description of the compactifications. A
canonical level structure will be sufficient for this \cite{A2}. Finally
the restriction to prime numbers $p$ in Theorem
0.2 is done to keep the technical difficulties to an acceptable level. I
believe that this restriction is not
essential. This result also supports a conjecture made in \cite{H}
for principal polarizations. This is
in so far interesting as the case
treated here cannot, as it could be in the case $p=1$,
be easily derived from known
results on ${\cal M}_2$.\\[3mm]

\noindent {\bf \large Acknowledgement: }  I am grateful to M.~Friedland and
G.K.~Sankaran for useful discussions.

\section{General type}
In this section we want to prove
\begin{theorem}
If $(d,n)=1$ and $n\ge 4$ then ${\cal A}_{1,d}(n)$ is of general type.
\end{theorem}

We shall work with the Voronoi compactification ${\cal A}^{\ast}_{1,d}(n)$
of ${\cal A}_{1,d}(n)$. Before we can prove
the theorem we need to know something about the coordinates of
${\cal A}^{\ast}_{1,d}(n)$ near a cusp. Recall that the codimension 1 cusps
are given by lines $l\subset \QQ^4$ up to the
action of the group
$\Gamma_{1,d}(n)$ and that the codimension 2 cusps are similarly given by
isotropic planes $h\subset \QQ^4$. For any
such $l$, resp. $h$ and any group $\Gamma$ we denote the lattice part of
the stabilizer of $l$, resp. $h$ in $\Gamma$ by
$P'_{\Gamma}(l)$, resp. $P'_{\Gamma}(h)$. These lattices have rank 1, resp. 3.

\begin{lemma}\label{lem1.2}
\begin{enumerate}
\item[{\rm(i)}]
For every line $l\subset\QQ^4$ there is an inclusion
$P'_{\operatorname{Sp}(4,\ZZ)}(l)\subset P'_{\Gamma_{1,d}}(l)$
with cokernel $\ZZ/d_1$ where $d_1|d$.
\item[{\rm(ii)}]
For every isotropic plane $h\subset\QQ^4$ there is an inclusion
$P'_{\operatorname{Sp}(4,\ZZ)}(h)\subset
P'_{\Gamma_{1,d}}(h)$ with cokernel $\ZZ/d_1\times\ZZ/d_2\times\ZZ/d_3$
where $d_i|d$ for $i=1,2,3$.
\end{enumerate}
\end{lemma}

\begin{Proof}
(i)\
Let $l$ be a line which corresponds to a given cusp of $\Gamma_{1,d}$. By
\cite[Satz 2.1]{FS} we may assume that
$$
l=\begin{pmatrix}
A & 0\\ C & D
\end{pmatrix} l_0\quad, \quad \begin{pmatrix}
A & 0\\ C & D
\end{pmatrix}=:M \in \mbox{Sp}(4,\ZZ)
$$
where $l_0=e_3=(0,0,1,0)$. The group $Q'(l_0)=M^{-1} P'_{\Gamma_{1,d}}(l)
M$ is a rank 1 lattice which fixes $l_0$. We
want to compare this to the rank 1 lattice
$$
P'_{\operatorname{Sp}(4,\ZZ)}(l_0)=\left\{\begin{pmatrix}
1 & 0 & s & 0\\
0 & 1 & 0 & 0\\
0 & 0 & 1 & 0\\
0 & 0 & 0 & 1
\end{pmatrix}; s\in \ZZ\right\}\subset \mbox{Sp} (4,\ZZ).
$$
Recall from \cite[Proposition I.1.16]{HKW} that every element $g$ in
$\Gamma_{1,d}$ fulfills the following
congruences
$$
g-{\bf 1} \in
\begin{pmatrix}
\ZZ     & \ZZ &  \ZZ   & d\ZZ\\
d\ZZ     & \ZZ & d\ZZ   & d\ZZ\\
\ZZ     & \ZZ & \ZZ    & d\ZZ\\
\ZZ     & \frac 1d\ZZ  &\ZZ & \ZZ
\end{pmatrix}
$$
Hence
$$
\begin{pmatrix}
^tD & 0\\
-^tC & ^tA
\end{pmatrix} g
\begin{pmatrix}
A & 0 \\
C & D\end{pmatrix}
=\begin{pmatrix}
\ast & S\\
\ast & \ast
\end{pmatrix}
,\quad S\in\begin{pmatrix}
\ZZ & \ZZ\\
\ZZ & \ZZ
\end{pmatrix}.
$$
In particular $P'_{\operatorname{Sp}(4,\ZZ)}(l_0)$ is contained in
$Q'(l_0)$ . Hence
$P'_{\Gamma_{1,d}}(l)/P'_{\operatorname{Sp}(4,\ZZ)}(l)\cong\ZZ/d_1$ for
some $d_1$. The claim $d_1|d$ follows since
$$
M
\begin{pmatrix}
1 & 0 & d & 0\\
0 & 1 & 0 & 0\\
0 & 0 & 1 & 0\\
0 & 0 & 0 & 1
\end{pmatrix}
M^{-1}\in P'_{\Gamma_{1,d}}(l).
$$

\noindent(ii)\
Again we can choose an element $M\in \operatorname{Sp}(4,\ZZ)$ such that
$h=M(h_0)$ where $h_0=e_3\wedge e_4$. Then
$$
Q'(h_0)=M^{-1} P'_{\Gamma_{1,d}}(h)M
$$
consists of elements of the form
$$
\begin{pmatrix}
1 & 0 & d_1\ZZ & d_2\ZZ\\
0 & 1 & d_2\ZZ & d_3\ZZ\\
0 & 0 & 1      & 0\\
0 & 0 & 0      & 1
\end{pmatrix}.
$$
By (i) we can conclude that $d_1, d_3\in\NN$. We claim that also
$d_2\in\NN$. To prove this recall that there is a
sublattice $L_0\subset L=\ZZ^4$ with $L/L_0\cong\ZZ/d$ such that
$\Gamma_{1,d}(L_0)\subset L$. (This is simply the
lattice spanned by $e_1, e_2, e_3, de_4$). Hence the same statement must be
true for $M^{-1} P'_{\Gamma_{1,d}}(h)M$,
but this implies that $d_2\in\NN$. The assertions $d_1|d$ and $d_3|d$
follow from (i) and $d_2|d$ follows again since
$$
M\begin{pmatrix}
1 & 0 & 0 & d\\
0 & 1 & d & 0\\
0 & 0 & 1 & 0\\
0 & 0 & 0 & 1
\end{pmatrix}M^{-1}\in P'_{\Gamma_{1,d}}(h).
$$
\hfill\end{Proof}

\noindent{\em Proof of the theorem.}
We consider the following maps of moduli spaces
$$
\diagram
& {\cal A}^{\operatorname{lev}}_{1,d}\dlto \drto & & {\cal A}_{1,d}(n)\dlto\\
{\cal A}_{1,1} & & {\cal A}_{1,d} &
\enddiagram
$$
The map ${\cal A}^{\operatorname{lev}}_{1,d}\rightarrow {\cal A}_{1,1}$
comes from the inclusion
$\Gamma^{\operatorname{lev}}_{1,d}\subset \operatorname {Sp}(4,\ZZ)$. (The
argument given in
\cite[Proposition I.1.20]{HKW} for $d$ prime goes through unchanged for all
$d$.). Note that
$\Gamma_{1,d}^{\operatorname{lev}}$ is not normal in
$\mbox{Sp} (4,\ZZ)$ and hence ${\cal
A}_{1,d}^{\operatorname{lev}}\rightarrow{\cal A}_{1,1}$ is not Galois. The
other
maps ${\cal A}_{1,d}^{\operatorname{lev}}\rightarrow {\cal A}_{1,d}$ and
${\cal A}_{1,d}(n)\rightarrow {\cal A}_{1,d}$
are Galois covers.\\
An essential ingredient in the proof is Igusa's modular form
$$
\Delta_{10}=\prod_{m\operatorname{ even}}\Theta^2_m(\tau)
$$
given by the product of the squares of all even theta null values. This is
a cusp form of weight 10 with respect to
$\mbox{Sp}(4,\ZZ)$. In fact it is, up to scalar, the unique weight 10 cusp
form with respect to $\mbox{Sp}(4,\ZZ)$.
Recall that it vanishes exactly along the $\mbox{Sp}(4,\ZZ)$-translates of
the diagonal
$$
\HH_1\times\HH_1=\left\{
\begin{pmatrix}
\tau_1 & 0\\
0 & \tau_3
\end{pmatrix};
\mbox{Im } \tau_1,\mbox{ Im } \tau_3>0
\right\}
\subset\HH^2
$$
where it vanishes of order 2. Since $\Gamma_{1,d}^{\operatorname{lev}}$ is
a subgroup of $\mbox{Sp}(4,\ZZ)$  we can
also consider $\Delta_{10}$ as a cusp form with respect to
$\Gamma_{1,d}^{\operatorname{lev}}$. Recall that for any
modular form $G$ and a matrix $M$ the slash-operator is defined by
$$
G|_k M:=\det(C\tau+D)^{-k}\  G(M \tau)\quad \left(M=\begin{pmatrix}A & B\\
C & D\end{pmatrix}\right).
$$
We consider the multiplicative symmetrization
$$
F_0:=\prod_{M\in\mbox{PSL}(2,\ZZ/d)}
\Delta_{10}|_{10}  M.
$$
It is straightforward to check that $F_0$ is a cusp form with respect to
$\Gamma_{1,d}$ of weight $10\mu(d)$ where
$$
\mu(d)=\frac 12|\mbox{SL}(2,\ZZ/d)|=\frac 12
d^3\prod_{p|d}(1-\frac{1}{p^2})\quad (d\ge 3),
$$
resp. $\mu(2)=6$. Clearly we can also consider $F_0$ as a cusp form with
respect to the smaller group $\Gamma_{1,d}(n)$.
Let $L$ be the ($\QQ$-)line bundle of modular forms of weight 1. By abuse
of notation we shall use the same notation
for whatever moduli space we are considering.\\

\noindent{\bf Claim 1:}
For every point $P$ on the boundary of ${\cal A}^{\ast}_{1,d}(n)$ the
modular form $F_0$ defines an element in
$m^{n\mu(d)}_P L^{10\mu(d)}$.

For points on the codimension 1 cusps this follows immediately from Lemma
1.2 (i) and $(n,d)=1$. To prove it in general
we consider an isotropic plane $h$ and the lattices
$N:=P'_{\Gamma_{1,d}}(h)$ and $N':=P'_{\Gamma_{1,d}(n)}(h)$. Let
$\sigma\in\Sigma_{\operatorname{vor}}$ be a 3-dimensional cone and let
$T_{\sigma}(N)$, resp. $T_{\sigma}(N')$ be the
corresponding affine parts in the toric variety
$T_{\Sigma_{\operatorname{vor}}}(N)$, resp.
$T_{\Sigma_{\operatorname{vor}}}(N')$. We claim that $\Delta_{10}$ defines
a function on the closure of the image of
$P'_{\Gamma_{1,d}}(h)\backslash \HH_2$ in $T_{\sigma}(N)$. First of all
$\Delta_{10}$ is a function on
$P'_{\Gamma_{1,d}}\backslash \HH_2$ by Lemma \ref{lem1.2}
(ii). Since $T_{\sigma}(N)$ is normal, it is enough to show that this
function extends to the codimension 1 boundary components. This follows
from Lemma \ref{lem1.2}
(i). Since $\Delta_{10}$ is a cusp
form it follows that $\Delta_{10}\in m_P$ for every point $P$ on the
boundary. By construction of $
F_0$ this gives
claim 1 in the case $n=1$. Since $(d,n)=1$ we have $N'=nN$ and this gives
the claim for general $n$.

Let ${\cal A}^{\ast}_{1,d}(n)$ be the Voronoi compactification of ${\cal
A}_{1,d}(n)$, i.e. the toroidal
compactification given by the second Voronoi decomposition of the cone of
semi-positive definite symmetric real
$(2\times 2)$
-matrices (in \cite{HKW} this was also called Legendre decomposition). If
$n\ge 3$ the group $\Gamma_{1,d}(n)$ is neat
(this follows from a general result of Serre which says that every
algebraic integer which is a unit and which is
congruent to 1 mod $n(n\ge 3)$ is equal to 1). In particular
${\cal A}_{1,d}(n)$ is smooth. The toroidal compactification
${\cal A}_{1,d}^{\ast}(n)$ will, in general, however have singularities.
These arise because the fan given by the Voronoi
decomposition is not always basic, i.e. there may be cones  which are not
spanned by elements of a basis of the
lattice. We can always choose a suitable subdivision of the fan given by
the Voronoi decomposition and in this way
construct a smooth resolution $\psi:\tilde{\cal A}_{1,d}(n)\rightarrow
{\cal A}_{1,d}^{\ast}(n)$ such that the
exceptional divisor is a normal crossing divisor.\\ Let
$\omega=d\tau_1\wedge d\tau_2\wedge d\tau_3$. It is well known
that, if $G$ is a weight $3 k$ cusp form which vanishes of order $\ge k$
along all 1-codimensional cusps, then
$G\omega^k$
 defines a $k$-fold canonical form on the smooth part of ${\cal
A}_{1,d}^{\ast}(n)$.\\

\noindent{\bf Claim 2:} The space of $k$-fold canonical forms which extends
to the smooth part of ${\cal
A}_{1,d}^{\ast}(n)$ grows (at least for sufficiently divisible $k$) as $c
k^3$ for some postive constant $c$.

To prove this claim recall that $K=3L-D$ on the smooth part of ${\cal
A}_{1,d}^{\ast}(n)$, where $L$ is the ($\QQ$-)
line bundle of modular forms of weight 1 and $D$ is the boundary, i.e. the
union of all 1-codimensional cusps. By claim
1 the form $F_0$ gives the equality
$$
10\mu(d)L=n\mu(d)D+D_{\operatorname{eff}}
$$
for some effective divisor $
D_{\operatorname{eff}}$ on the smooth part of ${\cal A}^{\ast}_{1,d}(n)$.
>From this we obtain
$$
-D=-\frac{10}{n} L+\frac{1}{n\mu(d)} D_{\operatorname{eff}}.
$$
Combining this equality with the expression for $K$ gives us
$$
K=\left(
3-\frac{10}{n}
\right)
L+\frac{1}{n\mu(d)} D_{\operatorname{eff}}.
$$
For $n\ge 4$ the factor in front of $L$ is positive and the claim follows
since $h^0(L^k)$ grows as $c k^3$.\\

\noindent{\bf Claim 3:} If $F_{3k}\omega^k$ defines a $k$-fold canonical
form on the smooth part of ${\cal
A}_{1,d}^{\ast}(n)$ then $\left(F_0^3\omega^{10\mu(d)}\right)^k\left(F_{3k}
\omega^k\right)$ extends to $\tilde{\cal
A}_{1,d}(n)$.
We
first notice that, since $n\ge 4$, the form $F^3_0\omega^{10\mu(d)}$
extends to the smooth part of ${\cal
A}^{\ast}_{1,d}(n)$. The following part of the argument follows closely the
proof of
\cite[Theorem 6.3]{S}. It is enough to prove that the form in question
extends to the generic point of each component of
the exceptional divisor. Let
$E$ be a component of the exceptional divisor of the resolution
$\tilde{\cal A}_{1,d}(n)\rightarrow {\cal
A}_{1,d}^{\ast}(n)$. It is enough to consider points which lie on only one
boundary component. We can choose local
analytic coordinates $z_1, z_2, z_3$ on an open set
$U$ such that
$E=\{z_1=0\}$. Recall that
$U$ is an open set in some toroidal variety $T_{\tilde{\Sigma}}(N')$ where
$N'=P'_{\Gamma_{1,d}(n)}(h)$ for some
isotropic plane
$h$ and
$\tilde{\Sigma}$ is a refinement of the fan
$\Sigma_{\operatorname{Vor}}$ defined by the Voronoi decomposition.
Moreover the coordinates of the torus are of the
form $t_i=e^{2\pi i a_i\tau_i}$ for some rational numbers $a_i$. A local
equation for $E$ is given by
$t_1^{b_1} t_2^{b_2} t_3^{b_3}$ for suitable $b_i$ and hence we can set
$z_1=t_1^{b_1} t_2^{b_2} t_3^{b_3}$. Since
${\partial z_1}/{\partial\tau_j}=2\pi i a_j b_j z_1$ we can conclude that
the order of
$J=\det\left({\partial \tau_i}/{\partial z_j}\right)$ along $E$ is
$v_E(J)\ge-1$. It follows again from claim 1 that
$$
v_E\left(F_0^3 J^{10\mu(d)}\right) \ge (3n-10)\mu(d).
$$
Therefore $\left(F_0^3\omega^{10\mu(d)}\right)$ defines a section of
$\mu(d)(10 K-(3n-10)E)$. By assumption
$F_{3k}\omega^k$ defines a section of $\psi^{\ast}\left(k K_{{\cal
A}^{\ast}_{1,d}(n)}\right)=k(K-\alpha E)$ where
$\alpha$ is the discrepancy of $E$. Altogether
$\left(F^3_0\omega^{10\mu(d)}\right)^k\left(F_{3k}\omega^k\right)$
defines a section of
$$
\begin{array}{r}
(10\mu(d)k K - k\mu(d)(3n-10)E) + k (K-\alpha E)=\\
=k\left[(10\mu(d)+1\right)K-\left(\mu(d)(3n-10)+\alpha)
E\right].
\end{array}
$$
All singularities here are cyclic quotient singularities. This follows from
Lemma \ref{lem1.2} and the fact that
$T_{\Sigma_{\operatorname{vor}}} (P'_{\operatorname{Sp}(4,\ZZ)}(h))$ is
smooth. Hence the singularities are
log-terminal, i.e. $\alpha > -1$. This implies that
$\mu(d)(3n-10)+\alpha > 0$ for $n\ge 4$ and thus the claim follows.

The theorem now follows easily from by combining claim 2 and claim 3.

\hfill$\Box$

\section{Ample canonical bundle}
It is the aim of this section to prove the following

\begin{theorem}
Let $p$ be an odd prime number and assume that $(n,p)=1$. The Voronoi
compactification $({\cal
A}_{1,p}^{\operatorname{lev}}(n))^{\ast}$ is smooth and has ample canonical
bundle if and only if $n\ge 5$.
\end{theorem}

We remark that this result is also known to be true for $p=1$ (cf.{\rm
\cite{Bo}, \cite{H}}). Before we prove this
theorem we recall the geometry of the spaces $({\cal
A}_{1,p}^{\operatorname{lev}})^{\ast}$ which was described in detail
in \cite{HKW}. The Tits building of the group
$\Gamma_{1,p}^{\operatorname{lev}}$ consists of
$1+(p^2-1)/2$ lines and $p+1$ isotropic planes. The lines consist of one
so-called central line and $(p^2-1)/2$
peripheral lines. If $D(l_0)$ is the closed boundary surface which belongs
to the central line, then there is a map
$K(p)\rightarrow D(l_0)$ which is an immersion, but not an embedding if
$p>3$. Here $K(p)$ is the Kummer modular
surface of level $p$, i.e. the quotient of Shioda's modular surface $S(p)$
by the involution which acts by $x\mapsto
-x$ on every fibre. For each peripheral boundary component $D(l)$ there
exists an isomorphism $K(1)\cong D(l)$ where
$K(1)$ is the Kummer modular surface of level 1.

If we add a level-$n$ structure clearly the number of inequivalent cusps
will increase. We shall, however, still speak
about central or peripheral cusps with respect to
$\Gamma_{1,p}^{\operatorname{lev}}(n)$ depending on whether this
defines a central or peripheral cusp with respect to
$\Gamma_{1,p}^{\operatorname{lev}}$. Now assume $n\ge 3$. Then
one shows exactly as in the proof of \cite[Theorem I.3.151]{HKW} that there
are immersions $S(np)\rightarrow D(l_c)$
if $l_c$ is a central cusp, resp. $S(n)\rightarrow D(l_p)$ if $l_p$ is a
peripheral cusp. The reason why we have
Shioda modular surfaces here instead of Kummer modular surfaces is that for
$n\ge 3$ the matrix -$\bf 1$ is not
contained in $\Gamma_{1,p}^{\operatorname{lev}}(n)$. It will be immaterial
for us whether these maps are immersions or
embeddings.

We shall write the boundary as
$$
D=\sum_{i\in I}D^i_c + \sum_{j\in J} D^j_p
$$
where $D^i_c$ are the central and $D^j_p$ the peripheral boundary components.

We recall the following well known facts about Shioda modular surfaces. For
$k\ge 3$ the surface $S(k)\rightarrow X(k)$
is the universal elliptic curve with a level-$k$ structure. The base curve
$X(k)$ is the modular curve of level $k$. It
has
$t(k)=\frac 12 k^2 \prod\limits_{p|k} (1-\frac 1{p^2})$ cusps and the fibre
of $S(k)$ over the cusps are singular of
type $I_k$, i.e. a $k$-gon of $(-2)$-curves. The genus of $X(k)$ equals
$1+(k-6)t(k)/12$ and the line bundle $L_{X(k)}$
of modular forms of weight 1 has degree $kt(k)/12$. The elliptic fibration
$\pi:S(k)\rightarrow X(k)$ has $k^2$
sections $L_{ij}$ which form a group $\ZZ/k\times \ZZ/k$. By $F$ we denote
a general fibre of $S(k)$. It is well known
(cf\cite[pp.78-80]{BH}) that
$$
\begin{array}{lrl}
K_{S(k)} &\equiv & \displaystyle\frac{k-4}{4} t(k)F,\\[2mm]
L^2_{ij} &=      & \displaystyle\frac {k}{12} t(k)=-\deg L_{X(k)}.
\end{array}
$$
Let $f_c:S(np)\rightarrow D_c$, resp. $f_p:S(n)\rightarrow D_p$ be the map
from $S(np)$ to a central, resp. from $S(n)$
to a peripheral boundary component. Since these maps are immersions we can
consider the normal bundles $N_c$, resp. $N_p$
of these maps.
\begin{proposition}
\begin{enumerate}
\item[{\rm(i)}]
$N_c\equiv -\displaystyle\frac{2}{np} L_{X(np)}-\frac{2}{np}
\sum\limits_{i,j\in\ZZ/{np}\times\ZZ/{np}}L_{ij}$
\item[{\rm(ii)}]
$N_p\equiv -\displaystyle\frac{2}n L_{X(n)}-\displaystyle\frac 2n
\sum\limits_{i,j\in\ZZ/n\times\ZZ/n} L_{ij}.$
\end{enumerate}
\end{proposition}

\begin{Proof}
We shall give the proof for the central boundary components and indicate
how it has to be adopted to the peripheral
boundary components. There is a natural action of the group $\sum
L_{ij}=\ZZ/{np}\times\ZZ/{np}$ on $S(np)$. It is an
easy calculation to check that this action is induced by elements of
$\Gamma_{1,p}^{\operatorname{lev}}(n)/\Gamma_{1,p}$. It follows that $N_c$
is invariant under the group
$\ZZ/{np}\times\ZZ/{np}$ and hence
$$
N_c\equiv aF+b\sum L_{ij}
$$
for some $a, b \in \QQ$ (cf.\cite{BH}). To determine $a,b$ we have to
compute the degree of the normal bundle $N_c$ on a
general fibre of $S(np)$ and on a section, e.g. the zero section $L_{00}$.\\
As a representative for a central cusp we can take the line $l_0=(0, 0, 1,
0)$. Assume $(n,p)=1$. We set
$$
\begin{array}{r}
\Gamma'_1(np):=\left\{
\begin{pmatrix} a & b\\ c & d\end{pmatrix}\right.
\in \operatorname{SL}(2,\ZZ); a,d  \equiv  1\mbox{ mod } np, c\equiv 0\mbox{
mod } n,\\
{\rule[-3mm]{0mm}{6mm}  \left.\begin{array}{cc}  \ & \ \\ \ & \ \end{array}
b  \equiv  0\mbox{ mod } np^2
 \right\} }.
\end{array}
$$

Note that by conjugation with $E=\mbox{diag} (1,p)$ the group
$\Gamma'_1(np)$ is conjugate to the principal subgroup
$\Gamma_1(np)$. Then by \cite[Proposition I.3.98]{HKW} the stabilizer
subgroup $P(l_0)$ of
$\Gamma^{\operatorname{lev}}_{1,p}(n)$ is given by
$$
P(l_0)=\left\{
\begin{pmatrix}
1 & k & s & m\\
0 & a & m & b\\
0 & 0 & 1 & 0\\
0 & c &-k & d
\end{pmatrix};\begin{pmatrix} a & b\\ c & d\end{pmatrix}\in\Gamma'_1(np),
k, s\in
n\ZZ, m\in pn\ZZ\}\right\}.
$$
The action of
$$
\begin{pmatrix}
1 & 0 & s & 0\\
0 & 1 & 0 & 0\\
0 & 0 & 1 & 0\\
0 & 0 & 0 & 1
\end{pmatrix}:\begin{pmatrix} \tau_1 & \tau_2\\ \tau_2 &
\tau_3\end{pmatrix}\mapsto\begin{pmatrix} \tau_1+s &
\tau_2\\ \tau_2
& \tau_3\end{pmatrix}
$$
gives rise to the partial quotient
$$
\begin{array}{ccl}
\HH_2 & \rightarrow & \CC^{\ast}\times \CC \times \HH_1\\
\begin{pmatrix}
\tau_1 & \tau_2 \\ \tau_2 & \tau_3\end{pmatrix}
&\mapsto& (t_1=e^{2\pi i\tau_1/n}, \tau_2, \tau_3).
\end{array}
$$
The other elements of $P(l_0)$ act as follows:\\
$
\begin{pmatrix}
1 & 0 & 0 & 0\\
0 & a & 0 & b\\
0 & 0 & 1 & 0\\
0 & c & 0 & d
\end{pmatrix}:\begin{pmatrix} \tau_1 & \tau_2\\ \tau_2 &
\tau_3\end{pmatrix}\mapsto
\begin{pmatrix} \tau_1-\tau_2(c\tau_3+d)^{-1} c\tau_2 &
\ast\\
\tau_2(c\tau_3+d)^{-1}
& (a\tau_3+b)(c\tau_3+d)^{-1}
\end{pmatrix},
$\\[2mm]
$
\begin{pmatrix}
1 & k & 0 & m\\
0 & 1 & m & 0\\
0 & 0 & 1 & 0\\
0 & 0 &-k & 1
\end{pmatrix}:
\begin{pmatrix} \tau_1 & \tau_2\\ \tau_2 & \tau_3\end
{pmatrix}
\mapsto
\begin{pmatrix} \tau'_1 &
\ast\\
\tau_2+k\tau_3+m
& \tau_3
\end{pmatrix},
$\\
${\ }\hspace{4.2cm}{\tau'_1=\tau_1+k^2\tau_3+2k\tau_2+km.}$\\[2mm]
The group
$$
P''(l_0)=\left\{
\begin{pmatrix}
1 & k & m\\
0 & a & b\\
0 & c & d
\end{pmatrix};
\begin{pmatrix}
a & b\\ c & d
\end{pmatrix}
\in \Gamma'_1(np), k\in n\ZZ, m\in pn\ZZ\right\}
$$
defines an action on $\CC\times\HH_1$ by\\[2mm]
$\begin{pmatrix}
1 & k & m\\
0 & 1 & 0\\
0 & 0 & 1
\end{pmatrix}
:(\tau_2, \tau_3)\mapsto (\tau_2+k\tau_3+m,\tau_3)
$\\[2mm]
$
\begin{pmatrix}
1 & 0 & 0\\
0 & a & b\\
0 & c & d\end{pmatrix}:(\tau_2, \tau_3)\mapsto (\tau_2(c\tau_3+d)^{-1},
(a\tau_3+b)(c\tau_3+d)^{-1}).
$\\
[2mm]
Then $D^0(l_0)=P''(l_0)\backslash\{0\}\times \CC\times\HH_1$ is the open
boundary surface associated to $l_0$ and
conjugation with $E=\mbox{diag}(1,p)$ shows that $D^{0}(l^0)\cong
S^{0}(np)$, the open part of $S(np)$ which does not lie
over the cusps.\\
A local equation of $D^{0}(l_0)$ in $\CC\times \CC\times\HH_1$ is given by
$t_1=0$ and hence $t_1/t_1^2$ is a local
section of the conormal bundle. Under the action of the group $P(l_0)$ this
transforms as follows:
\begin{enumerate}
\item[(1)]\qquad $t_1/t_1^2 \mapsto t_1/t_1^2\  e^{2\pi i[-k^2
\tau_3-2k\tau_2]/n},$\\

\item[(2)]\qquad $t_1/t_1^2\mapsto t_1/t_1^2\  e^{2\pi
i\left(-\frac{c\tau_2^2}{c\tau_3+d}\right)/n}.$
\end{enumerate}
We can use the formulae (1) and (2) to determine the coefficients $a$ and
$b$. We first determine the degree of $N_c$
on a general fibre $F$. Since $k\in n\ZZ,\ m\in pn \ZZ$ the fibre of
$S(np)$ over the point $[\tau_3]\in X(np)$ is
given by $E_{[\tau_3]}=\CC/(\ZZ n\tau_3+\ZZ np)$. The standard theta
function $\vartheta(\tau_3, \tau_2)$ defines a line
bundle of degree $n^2p$ on $E_{[\tau_3]}$ and transforms as follows
\begin{enumerate}
\item[(3)] $\displaystyle\vartheta(\tau_3, \tau_2+k \tau_3+m)=e^{2\pi
i[-\frac 12 k^2\tau_3-k\tau_2]}.$
\end{enumerate}
Comparing formulae (1) and (3) we find that the degree of $N_c$ on $F$
equals $-2np$. Since  we have $n^2p^2$ sections
$L_{i j}$ it follows that $b=-2/np$. To determine the coefficient $a$ we
have to compute the degree of $N_c$ on the
zero section $L_{00}$. Since $L_{00}^2=-\deg L_{X(np)}$ we must show that
this degree is $0$. There are two ways to see
this. The first is to use formula (2) and specialise it to $\tau_2=0$. One
then has to show that this description extends
over the cusps which is an easy local calculation. Alternatively one can
proceed as follows: The section $L_{00}$ is the
transversal intersection of $D_c$ with the closure of the image of the
diagonal $\HH_1\times \HH_1 \subset \HH_2$ which
parametrizes products. This closure is isomorphic to $X(n)\times X(np)$ and
$L_{00}$ is equal to $\{\mbox{cusp}\}\times
X(np)$. Hence the normal bundle of $L_{00}$ in $X(n)\times X(np)$ is
trivial and by adjunction
$$
K_{L_{00}}=K|_{L_{00}}+L_{00}|L_{00}
$$
where $K$ is the canonical bundle of $({\cal
A}^{\operatorname{lev}}_{1,p}(n))^{\ast}$. Using the fact that $K=3L-D$
and pulling this back to $S(np)$ we obtain
$$
K_{L_{00}}=\left( 3L_{X(np)}-t(np)F-N_{c}+ L_{00}\right)|_{L_{00}}.
$$
Since $\deg K_{L_{00}}=t(np)(np-6)/6$ a straightforward calculation shows
that the degree of $N_c|_{L_{00}}$ is equal
to $0$.

The calculation for $N_p$ is essentially the same. The only differences are
that $t_3=e^{2\pi i\tau_3/np^2}$ and that
the fibre of $S(n)$ over $[\tau_1]\in X(n)$ is equal to
$E_{[\tau,_3]}=\CC/(\ZZ np\tau_1+np\ZZ).$
\hfill$\Box$\\

\noindent{\em Proof of the theorem:} We first remark that $\left({\cal
A}^{\operatorname{lev}}_{1,p}(n)\right)^{\ast}$ is
smooth under the assumptions made. Since $n\ge 4$ is the group
$\Gamma^{\operatorname{lev}}_{1,p}(n)$ is neat.
Therefore it is enough to show that for a given cusp $h$ the toroidal variety
$T_{\Sigma_{\operatorname{vor}}}(P'_{\Gamma_{1,p}^{\operatorname{lev}}(n)}(h))$
 is smooth. If $n$ and $p$ are coprime, the lattice
$P'_{\Gamma_{1,p}^{\operatorname{lev}}(n)}(h)$ is simply
$n$ times the corresponding lattice in $\Gamma^{\operatorname{lev}}_{1,p}$.
Hence every
$\sigma\in\Sigma_{\operatorname{vor}}$ is spanned over $\RR$ by a basis of
the lattice and this implies that
$T_{\sigma}$ is smooth.
\\The next observation is that the condition $n\ge 5$ is necessary. We have
already remarked that the closure of the
diagonal $\HH_1\times \HH_1\subset \HH_2$ parametrizing split abelian
surfaces is isomorphic to $X(n)\times X(np)$.
Consider a curve $C=X(n)\times\operatorname{\{point\}}$. Then
$$
K|_C=(3L-D)|_C=3L_{X(n)}-X_{\infty}(n)
$$
where $X_{\infty}(n)$ is the divisor of cusps on $X(n)$. Hence
$$
K.C=\frac n4 t(n)-t(n)
$$
and this is positive if and only if $n\ge 5$.\\
We shall now assume $n\ge 5$. Let $C$ be an irreducible curve which is not
entirely contained in the boundary, i.e.
$C\cap{\cal A}^{\operatorname{lev}}_{1,p}(n)\neq\emptyset$ and consider a
point $[\tau]\in C$. Choose some
$\varepsilon>0$ with $\varepsilon<3n/5$. By Weissauer's result \cite[p.
220]{We} we can find a cusp form $F$ with respect
to
$\operatorname{Sp}(4,\ZZ)$ such that $F(\tau)\neq 0$ and $o(F)\ge
1/(12+\varepsilon)$. Here $o(F)$ is the order of
$F$, i.e. the vanishing order of $F$ divided by the weight of $F$. Let $m$
and $k$ be the vanishing order, resp. the
weight of $F$. Since $\Gamma^{\operatorname{lev}}_{1,p}(n)\subset
\operatorname{Sp}(4,\ZZ)$ the form $F$ is also a
modular form with respect to $\Gamma^{\operatorname{lev}}_{1,p}(n)$. In
terms of divisors this gives us
$$
kL=mnD+D_{\operatorname{eff}},\quad C\not\subset
D_{\operatorname{eff}}
.$$
Here $D_{\operatorname{eff}}$ contains in particular multiples of
peripheral boundary components since the vanishing
order of $F$ along these boundary components is at least $np^2$. From the
above formula we find that
$$
\left(\frac k{mn} L-D\right). C=\frac 1{mn} D_{\operatorname{eff}}. C\ge 0.
$$
Since $L.C>0$ we find that $K.C>0$ provided $3>k/mn$. But this follows
immediately from the inequalities $m/k\ge
1/(12+\varepsilon)$ and $\varepsilon<3n/5$.\\
It remains to prove that the restriction of $K$ to every boundary component
is ample. Let $D_0$ be a boundary component
and set $D'_0=D-D_0$. We have already observed that there is an immersion
$f:\tilde{D}\rightarrow D_0$ which is the
normalization. The surface $\tilde{D}$ is either isomorphic to $S(np)$ or
to $S(n)$ depending on whether we have a
central or a peripheral boundary component. The map $f$ embeds every
component of a singular fibre. The image of such a
component in
$\left({\cal A}^{\operatorname{lev}}_{1,p}(n)\right)^{\ast}$ is a $\PP^1$.
Away from $\{0,\infty\}$ this line is either
the intersection of 2 different boundary components or 2 branches of $D_0$
intersecting transversally. In either case we
have the
$$
f^{\ast}K=f^{\ast}(3L-D'_0-D_0)=3L_X-F_{\infty}-N_f.
$$
Here $L_X$ is either $L_{X(np)}$ or $L_{X(n)}$ depending on the type of the
boundary component, the divisor $F_{\infty}$
is the union of the singular fibres and $N_f$ is the normal bundle of the
immersion $f$. Let $k=np$ or $n$. Then
$$
\mbox{deg}(3L_X-F_{\infty})=\frac 14 kt(k)-t(k)>0
$$
for $n>4$. Hence $(3L_X-F_{\infty}). C\ge0$ for every curve $C$ and
$(3L_X-F_{\infty}). C>0$ unless $C$ is contained in
a union of fibres. It follows immediately from our proposition that
$-N_f.C>0$ for every curve $C$ which does not contain
a section
$L_{ij}$. Since $-N_f.L_{ij}=0$ we can conclude that $f^{\ast}K.C>0$ for
every curve $C$.
\hfill
\end{Proof}
\begin{remark}
The above proof can also be adapted to show that $K$ is nef for $n=4$. We
had already seen that $K$ is not ample in this
case. In other words
$\left({\cal A}^{\operatorname{lev}}_{1,p}(4)\right)^{\ast}$ is a minimal,
but not a canonical model for $p=1$ or $p\ge
3$ prime.
\end{remark}


%
\bibliographystyle{alpha}

\vspace{1cm}
Klaus Hulek\\
Universit\"at Hannover
\\Institut f\"ur Mathematik\\
Postfach 6009\\
D-30060 Hannover\\
Germany\\
hulek@math.uni-hannover.de

\end{document}